\documentclass[11pt]{amsart}
\usepackage{amsmath, amsthm, amsfonts, amssymb, graphicx}
\usepackage{setspace}
\usepackage{fullpage}
\usepackage{amscd}
\usepackage{amsaddr}
\theoremstyle{plain} \numberwithin{equation}{section}

\newtheorem*{theorem*}{Theorem}

\def\R{\mathbb R}
\def\00{(0,0)}
\def\11{(1,1)}

\begin{document}

\title{Gradient descent in higher codimension}
\author{Y.\ Cooper}
\email{yaim@math.ias.edu}
\date{\today}
\maketitle

\vspace{-.1in}  

\begin{abstract}
We consider the behavior of gradient flow and of discrete and noisy gradient descent.  It is commonly noted that the addition of noise to the process of discrete gradient descent can affect the trajectory of gradient descent.  In \cite{codim1}, we observed such effects.  There, we considered the case where the minima had codimension 1.  In this note, we do some computer experiments and observe the behavior of noisy gradient descent in the more complex setting of minima of higher codimension.  
\end{abstract}

\section{Introduction}

In this note, we explore the behavior of gradient flow in a landscape that contains minima of different geometries, as well as the behavior of several kinds of noisy discrete gradient descent the same setting.  We observe that the addition of noise to the process can affect the trajectory of gradient descent.

We will consider gradient descent on the function 
$$f(x,y) = \sin(\pi x) \sin(2 \pi x)  \cos(\pi y) \cos( 2 \pi y).$$  

\begin{figure}[h]
\includegraphics[width=3.5in]{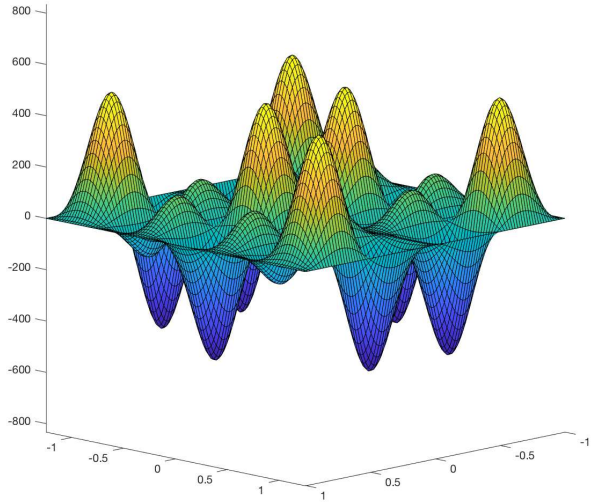}
\caption{We consider gradient descent on the function $f(x,y) = \sin(\pi x) \sin(2 \pi x)  \cos(\pi y) \cos( 2 \pi y)$, which has wells of two different depths.
\label{graph}}
\end{figure}

This function is doubly periodic, and has two kinds of wells: a deeper wider well, and a shallower narrower well.  This function has more features than the ones considered in \cite{codim1}.  For $f(x,y)$ there are many saddle points where gradient descent can slow down near.  This is possible also in one-dimensional examples, but was not an issue in the cases considered in \cite{codim1}.  For $f(x,y)$ the arrangement of the deep and shallow wells is not in a checkerboard pattern, but a more complex one, where they come in pairs of deep and shallow wells, as shown in Figure \ref{graph}.  

In general, the behavior of gradient descent on high dimensional landscapes is likely more complex than the behavior of gradient descent in low dimensional landscapes, in part because the landscapes themselves can become much more complex.  It is hard for us to visualize graphs of functions from $\R^n$ to $\R$ when $n$ becomes larger than 2.  So in this note, we focus on the case $n=2$, and try to understand ways in which the addition of noise can affect the trajectories of gradient descent. 

\subsection{Acknowledgements}

The author would like to thank Avi Wigderson for conversations that inspired this experiment, and Nathaniel Bottman for help with the experiments.

\section{Gradient descent on $f(x,y)$: Theory}
We will consider gradient descent on the function 
$$f(x,y) = \sin(\pi x) \sin(2 \pi x)  \cos(\pi y) \cos( 2 \pi y).$$  

This function is doubly periodic, and has two kinds of minima, a deeper minimum with depth approximately 0.77, and a shallower one of depth approximately 0.21.  The width of the deep wells is approximately 0.25, and the width of the shallow wells is approximately 0.125.

\subsection{Gradient flow}

We expect that if we start at a point $p_0$ and follow the trajectory under gradient flow, we will come to a critical point.  $f(x,y)$ has minima, maxima, and saddle points.  The trajectory will never flow to a maximum, but it can end near a saddle point or a minimum.  We will not consider the behavior around saddle points in this note.  Rather, we will study the relative tendency to end near a deeper or shallower minimum.  Given uniform initialization in the region of study, we expect the relative probabilities to be approximately equal to the relative areas of the basins of attraction.

In the lower dimensional case of \cite{codim1}, this was relatively easy to calculate.  In the case of $f(x,y)$, with minima of codimension 2, it is already difficult to calculate directly.  So we will approximate the relative areas by computer experiment running discrete gradient descent with no noise and small step size.  

\subsection{Discrete gradient descent with $\epsilon$-jitter}\label{with jitter}
We are interested in observing the behavior of discrete gradient descent with added noise.  Following \cite{codim1}, in this paper we will experimentally observe the behavior of discrete gradient descent with $\epsilon$-jitter, defined as follows.   

We begin at some initial position $p_0$.  At the $t+1^{st}$ step, we let
$$
p_{t+1} = p_t - \tau \nabla L (p_t) - (\epsilon_{t})
$$
where $\epsilon_{t}$ are drawn from a gaussian distribution with norm 0 and standard deviation $\epsilon$.  

\section{Discrete gradient descent on $f(x)$: Computer experiments}

We now turn to some computer experiments.  In all the experiments of this section, we run variants of gradient descent on the function $f(x,y)$ over the region $x \in [-1,1]$ and $y \in [-1.25,1.25]$ and initialize from the uniform distribution on that region.

\subsection{Gradient flow}
\begin{figure}[h]
\includegraphics[width=3.4in]{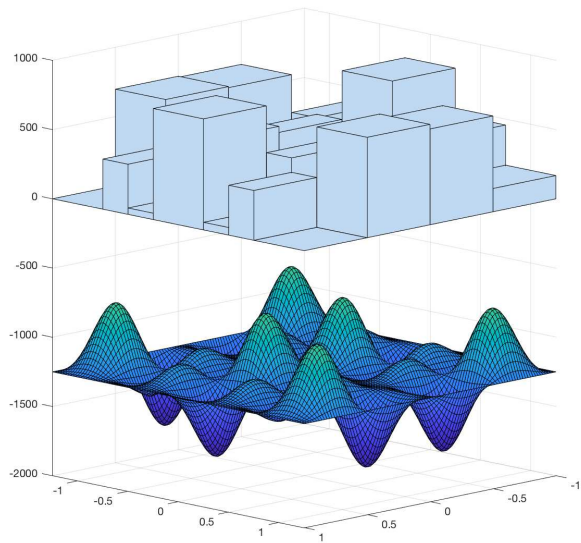}
\caption{In this experiment, we implement gradient descent with step size $\tau = 0.001$ on the function $f(x,y) = \sin(\pi x) \sin(2 \pi x)  \cos(\pi y) \cos( 2 \pi y)$.  We randomly initialize from the uniform distribution on this interval, and run the experiment 10,000 times.  Above is a histogram which shows the number of times the process ended in each well.}
\end{figure}

First, we approximate gradient flow by implementing discrete gradient descent with step size $\tau = 0.001$.  When we do this, we find that the trajectories usually end in one of the two types of minima, but also regularly get stuck near saddle points.  The ratio $r$ of the probability of ending in a shallow well to the probability of ending in a deep well is approximately $0.78$, and the probability of the trajectory ending within the test area is $0.89$.

\subsection{Discrete gradient descent with $\epsilon$-jitter}\label{with jitter}

\begin{figure}[h] \label{2}
\includegraphics[width=3.4in]{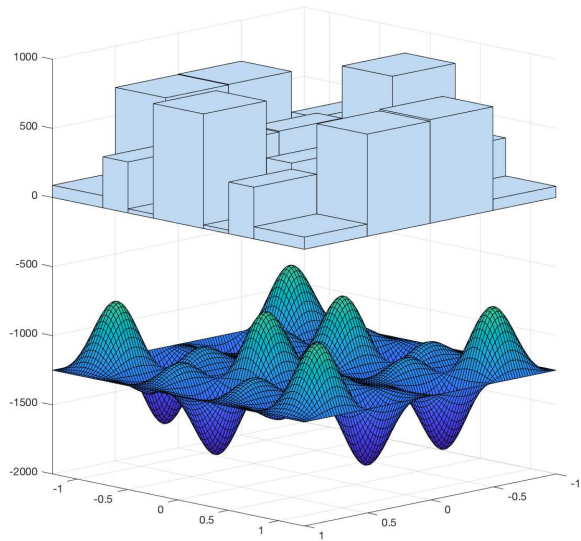} 
\caption{This histogram shows the number of times the process ended in each well with $\epsilon = 0.01$.  The ratio $r$of the probability of landing in a shallow well to the probability of landing in a deep well was 0.70.}
\end{figure}

Now, we turn to computer experiments where we run discrete gradient descent with $\epsilon$-jitter.  First, we fix the step size $\tau = 0.01$ and vary the amount of noise $\epsilon$ that we add.  In all the experiments in this subsection, we run 10,000 trials, and for each, we take 500 steps.  

To visualize the results, we show a histogram of the number of times the process ends in each region.  We have divided the area of study along the lines where $f(x,y)=0$.  So some of the resulting regions correspond to hills and some to wells.  There are also saddle points along the grid lines.  Below the histogram, we show the graph of $f(x,y)$.  We also note the ratio $r$ between the probability of ending in a shallow well to the probability of ending in a deep well, as well as the fraction $\phi$ of trials that stay within the region of study.

In our first experiment, we use $\epsilon = 0.01$.  We find that even adding a small amount of noise in this setting has an effect on the relative probabilities of ending in a deep or shallow well.  In particular, compared to the implementation of gradient descent with no noise, the noisy process is more likely to find the deeper minima.  The results are shown in Figure \ref{2}.  The ratio $r$ found was $0.71$, and the probability $\phi$ was $0.87$.

In the second run, we use $\epsilon = 0.05$.  At this level of noise, the effect becomes even more pronounced, and the relative probability of finding a shallow minimum rather than a deep one becomes $r=0.13$, while the probability of staying within the test area is $\phi = 0.81$.

\begin{figure}[h] \label{3}
\includegraphics[width=3.4in]{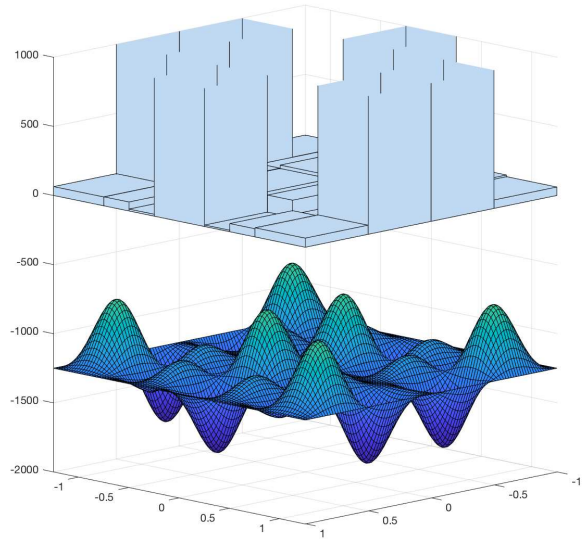}
\caption{This histogram shows the number of times the process ended in each well with $\epsilon = 0.05$.  The ratio $r$ of the probability of landing in a shallow well to the probability of landing in a deep well was 0.13.}
\end{figure}

Next, we take $\epsilon = 0.09$.  Now, the effect is so strong that the process almost always finds the deeper wells.  The relative probability $r$ drops to $0.03$.  However, the higher noise level also causes more of the trajectories to end outside of the test area, and $\phi$ decreases to $0.67$.

\begin{figure}[h] \label{4}
\includegraphics[width=3.4in]{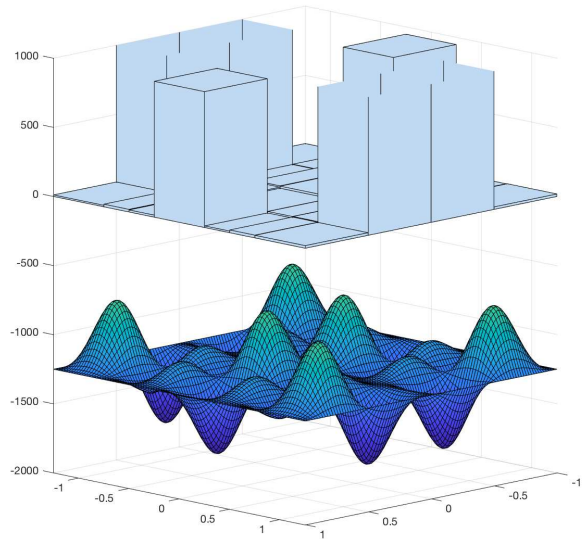}
\caption{This histogram shows the number of times the process ended in each well with $\epsilon = 0.09$.  The ratio $r$ of the probability of landing in a shallow well to the probability of landing in a deep well was 0.03.}
\end{figure}

Finally, we take $\epsilon = 0.15$.  This level of noise is sufficient that most of the trajectories end outside the test area, and $\phi$ drops to $0.21$.  The effect on biasing the process toward deep minima persists, but is less strong, and $r$ rises to $0.33$.  However, this level of noise is high enough that the process no longer behaves similarly to gradient flow, and rather is a process dominated by noise.

\begin{figure}[h] \label{5}
\includegraphics[width=3.4in]{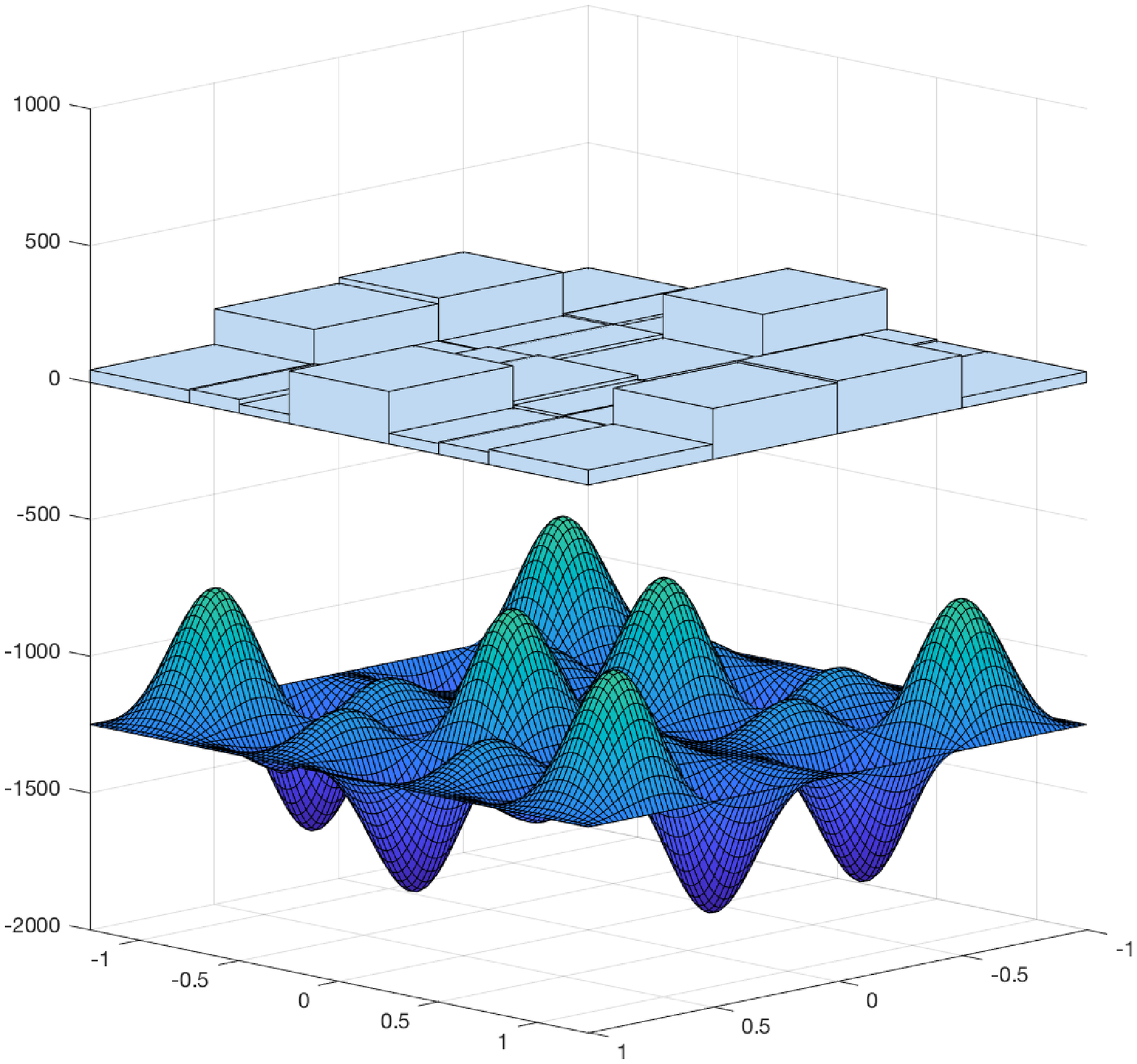}
\caption{This histogram shows the number of times the process ended in each well with $\epsilon = 0.15$.  The ratio $r$ of the probability of landing in a shallow well to the probability of landing in a deep well was 0.33.}
\end{figure}

To make a more systematic study of how the strength of this effect depends on the step size $\tau$ and the noise $\epsilon$, we now run a number of experiments at different values of $\tau$ and $\epsilon$, and observe how the behavior of noisy gradient descent varies.  We do this for step sizes $\tau = 0.001, 0.01, 0.02, 0.04, $ and $0.06$.  For each fixed $\tau$, we vary the noise $\epsilon$ from $0$ to $0.3$.  For each choice of $\tau$ and $\epsilon$, we do 500 trials, each consisting of a randomly initialized trajectory minimized by discrete gradient descent with $\epsilon$ jitter, with run for 500 steps with step size $\tau$.  

\begin{figure}[h] \label{6}
\includegraphics[width=6in]{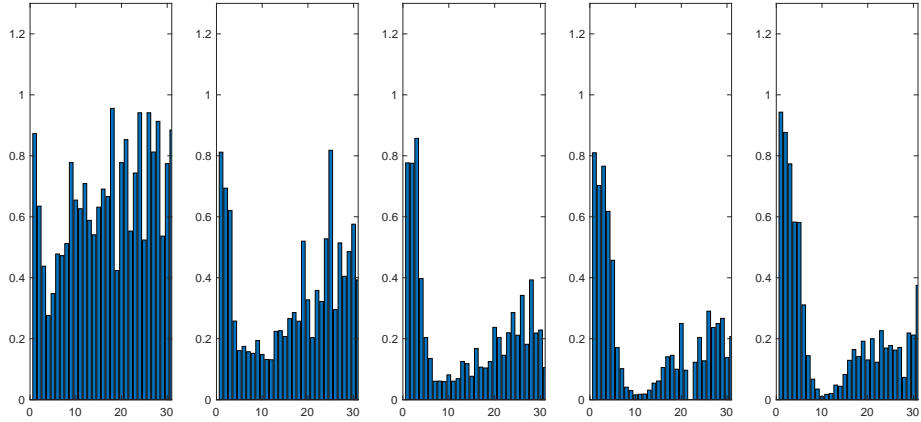}
\caption{These charts, from left to right, concern the behavior of noisy gradient descent on $f(x,y) = \sin(\pi x) \sin(2 \pi x)  \cos(\pi y) \cos( 2 \pi y)$, with step size $\tau = 0.001, 0.01, 0.02, 0.04$, and $0.06$.  Each chart shows, for that step size, how the ratio of the probability of ending in a deep vs a shallow well varies as the magnitude $\epsilon$ of the added noise varies from $0$ to $0.3$.}
\end{figure}

We display the results in Figure \ref{6}.  There are 5 charts displayed, one for each choice of step size $\tau$.  In each chart, the noise $\epsilon$ increases along the $x$-axis, while the $y$-axis displays the ratio $r$ between the probability of ending in a shallow well to the probability of ending in a deep well.  

We find that compared to the codimension 1 case studied in \cite{codim1}, in this setting the effect of adding noise is more robust, across a range of step sizes.  The smallest step size $\tau = 0.001$ is a discrete approximation of gradient flow, and in this setting, the majority (more than 80 \%) of the trajectories that began in the test area stayed in the test area.  In contrast, the largest step size $\tau = 0.06$ is large enough that in the codimension 1 case, even mild addition of noise caused the process to become fairly unstable and for many of the trajectories to leave the test area.

In this setting with minima of codimension 2, the addition of noise substantially biases the process toward deep minima, at all step sizes considered.  The effect is smallest for the smallest step size $\tau = 0.001$.  Both the width of the gap and the depth of the gap in the histogram is smallest in this case, and we observe that the effect is very sensitive to the exact amount of noise added.  

However, for all other step sizes considered here, the addition of noise reliably biases the process toward deep minima.  In the cases $\tau = 0.01$ to $0.06$, the strength of the effect is fairly robust to the exact amount of noise injected, and for $\tau = 0.04$ and $\tau =0.06$, at some values of $\epsilon$ the added noise causes the ratio $r$ of the probability of ending in a shallow well to the probability of ending in a deep well to drop almost to 0.  

\section{Discussion}

In \cite{codim1}, we compared the behavior of discrete gradient descent approximating gradient flow to the behavior of noisy discrete gradient descent in landscapes with minima of different geometries.  We observed that in that setting, noise can bias the procedure toward deeper and wider wells, as is often suggested.

In this note, we consider the same question in the case of a landscape in dimension 2, with minima of codimension 2.  We observed that in this setting as well, noise can bias the procedure toward deeper and wider wells.  There were some small differences between the experiments in dimension 1 and dimension 2.  For one, in the two dimensional example, the process sometimes ended near saddle points, which didn't happen in the one dimensional examples we considered in \cite{codim1}, though in principle it could have.  

We also noticed that in the experiments in this paper, the effect of noise seemed to be more robust to the exact amount of noise injected than it was in the one dimensional cases we tried.  We do not know if that is a property of working in a higher dimension, or if it was a property of the specific functions that we were experimenting with.  We leave that question to future work.

\end{document}